\documentclass{article}
\usepackage{ijcai20}

\usepackage{times}
\usepackage{soul}
\usepackage{url}
\usepackage[hidelinks]{hyperref}
\usepackage[utf8]{inputenc}
\usepackage[small]{caption}
\usepackage{graphicx}
\usepackage{amsmath}
\usepackage{amsthm}
\usepackage{booktabs}
\usepackage{algorithm}
\usepackage{algorithmic}
\usepackage{subfigure}

\title{Stochastic Item Descent Method\\for Large Scale Equal Circle Packing Problem}

 \author{
 Kun He$^1$
 \and
  Min Zhang$^1$\thanks{Corresponding author.} \and
 Jianrong Zhou$^{1}$\and
 Yan Jin$^1$\and
 Chu-min Li$^2$
 \affiliations
 $^1$School of Computer Science and Technology, Huazhong University of Science and Technology,
 China\\
 $^2$MIS, Universit\'e de Picardie Jules Verne, 
 France
 \emails
 \{ brooklet60, m\_zhang, yukihana0416,jinyan\}@hust.edu.cn,
 chu-min.li@u-picardie.fr
 }

\begin{document}

\maketitle

\begin{abstract}
 Stochastic gradient descent (SGD) is a powerful method for large-scale optimization problems in the area of machine learning, especially for a finite-sum formulation with numerous variables. In recent years, mini-batch SGD gains great success and has become a standard technique for training deep neural networks fed with big amount of data. Inspired by its success in deep learning, we apply the idea of SGD with batch selection of samples to a classic optimization problem in decision version. Given $n$ unit circles, the equal circle packing problem (ECPP) asks whether there exist a feasible packing that could put all the circles inside a circular container without overlapping. 
 Specifically, we propose a stochastic item descent method (SIDM) for ECPP in large scale, which randomly divides the unit circles into batches and runs Broyden-Fletcher-Goldfarb-Shanno (BFGS) algorithm on the corresponding batch function iteratively to speedup the calculation. We also increase the batch size during the batch iterations to gain higher quality solution. Comparing to the current best packing algorithms, SIDM greatly speeds up the calculation of optimization process and guarantees the solution quality for large scale instances with up to 1500 circle items, while the baseline algorithms usually handle about 300 circle items. 
 The results indicate the highly efficiency of SIDM for this classic optimization problem in large scale, and show potential for other large scale classic optimization problems in which gradient descent is used for optimization.

\end{abstract}

\vspace{-1em}
\section{Introduction}
\label{010intro}
Stochastic gradient descent (SGD) method~\cite{robbins1951} has gained great success in the area of machine learning~\cite{10.1007/978-3-7908-2604-3_16,doi:10.1137/16M1080173}. Especially for deep learning tasks, mini-batch SGD has become a standard technique for the training of deep neural networks fed with big amount of  data~\cite{DBLP:books/daglib/0040158,LeNet1998} . 
Inspired by its successful application for such big, complex optimization problems, in this work, we consider a classic global optimization problem well-studied in the area of operations research for over 30 years\cite{doi:10.1080/0025570X.1967.11975768}, and apply the idea of batch gradient descent (BGD) for this problem in large scale.  

Specifically, we consider the equal circle packing problem (ECPP) in decision version, 
the purpose of which is to answer whether a dense arrangement of $n$ unit circles without overlapping (i.e. feasible) in a circular container of fixed radius.  
If we already have an efficient algorithm for the decision version, the optimal version of minimizing the container radius for feasible packings can be solved efficiently by combining divide and conquer on the container radius. 
Our motivation is how to design an algorithm that is very fast so as to address the problem in large scale where hundreds and thousands of unit circles are considered.

Finding the optimal solution of ECPP with plenty number of circles is known to be NP hard, even the search of a suboptimal solution is still very challenging. Many researchers design heuristic algorithms to find a suboptimal packing pattern. In recent years, the quasi-physical energy based method was proposed which could solve ECPP in optimal version with up to a hundred items. Many quasi-physical researches regard each circle as an elastic item and treat the container as a rigid hollow container~\cite{He:2013:CQO:2560985.2561717,He2015,HE201826}. If two items, or an item with the container are squeezed against each other, the whole system would have elastic potential energy, and by gradient descent method like Broyden-Fletcher-Goldfarb-Shanno (BFGS) algorithm \cite{liu1989limited} we can reduce the potential energy of the system so as to remove the overlapping. Then some Basin-hopping strategy is used to jump out of the local optimal trap where overlaps still exist.

On the other hand, SGD is a classic first-order optimization algorithm widely used in large scale machine learning problems due to its low computational cost and modest precision~\cite{10.1007/978-3-7908-2604-3_16,Bottou08thetradeoffs}. In the training of deep neural networks, SGD plays a key role in the optimization process, and promotes the great success of deep learning. In the iteration of SGD, it randomly selects a sample and then optimizes the loss function corresponding to the current sample. 

Inspired by the success of SGD in deep learning, can we apply this idea to the classic optimization problem of ECPP? Specifically, can we randomly select a unit circle and optimize the corresponding optimization function? If each time a batch of circles are selected for gradient descent by fixing other circles, then we have a batched version of SGD. As quasi-Newton methods have been shown superior to first-order gradient descent method for various circle packing problems, we choose a quasi-Newton method like BFGS and combine it with random batch, and design a batched version of stochastic BFGS for ECPP. 

Therefore, we propose a novel approach called stochastic item descent method (SIDM), which can find dense layouts for large scale ECPP. SIDM accelerates the search process, especially for a large number of unit circles. In addition, after attaining a local minimum or saddle points, we improve the hopping strategy in the current best solution~\cite{HE201826}, in which we gradually increase the shrinking radius of container during the iteration to find better solutions. Comparing to state-of-the-art algorithms that can only address small-scale ECPP within reasonable time, SIDM can address up to $n = 1500$ instances and reach current best solution reported on the packomania website \footnote{http://www.packomania.com}. 

Our main contributions are listed as follows:
\vspace{-3pt}
\begin{itemize}
\item The proposed novel method SIDM can speed up the process of reaching the local minimum or saddle point, which is the main computation load of the ECPP.
\item We improve the basin-hopping procedure of the existing strategy used to escape  suboptimal layouts, and shrink the radius of the container more flexibly.
\item Experiments demonstrate that SIDM can greatly accelerate the computation while maintaining the state-of-art packing quality. 
\end{itemize}
\vspace{-3pt}


\vspace{-0.5em}
\section{Related Work}
\label{011rela}
In the literature, most researchers address the optimal version of ECPP
that requires to find the smallest container radius for all items. But they usually solve the decision version of ECPP as a sub-problem
and then use binary search (divide and conquer on the container radius) so as to find a possible smallest container radius for feasible packing. 
The efficiency and effectiveness of the overall algorithm mainly depend on the algorithm on the decision version.
Thus, in this work, we focus on improving the efficiency of the sub-algorithm for the decision version while maintaining the same effectiveness. 
And in the following, we provide an overview for the ECPP in optimal version. 

ECPP is a well studied problem since  1960's~\cite{doi:10.1002/mana.19690400110}. Mathematicians found the optimal packing pattern for $1 \leq n \leq 13$~\cite{doi:10.1002/mana.19690400110,Melissen199415,Fodor2000401,Fodor2003431} and $n = 19$~\cite{Fodor1999139}. 
However, it is very hard to mathematically find optimal solutions for bigger $n$, and mathematicians only found suboptimal packing patterns for $n \leq 25$~\cite{doi:10.1002/mana.19690400110,doi:10.1080/0025570X.1971.11976122,Reis1975}.

To earn a good trade-off between the computation efficiency and solution equality, greedy based heuristic algorithms performance well for $n\le100$. 
Graham \textit{et al.} proposed methods that simulate repulsion forces and billiards to iteratively search for global optimal layout \cite{Graham:1998:DPC:282687.282708}, and found suboptimal solution for $25\leq n \leq 65$. Akiyama \textit{et al.} obtained dense layout for $n = 70, 73, 75, 77, 78, 79, 80$ by a greedy algorithm \cite{Akiyama2003}. Then, Grosso \textit{et al.} proposed a monotonic basin hopping algorithm that improved many solutions for $66 \leq n \leq 100$~\cite{Grosso:2010:SPP:1773612.1773615}. 


For heuristic approaches, a typical way is to transform ECPP into a discrete optimization problem, i.e. putting the unit circles into the container one by one~\cite{chen2018greedy}, and then incorporating some search methods to improve the solution. Beam search algorithm~\cite{AKEB20091513} 
and greedy heuristic algorithm~\cite{chen2018greedy} have been proposed, which are all based on max hole degree method~\cite{Huang2003}. However, the solution quality is rather limited. 

Another approach is to formulate ECPP into a continuous optimization problem, that is, put all circles into the container allowing overlapping, use gradient based optimization algorithms to constantly adjust positions of the unit circles, and shrink the container radius for the next round of search if feasible solution is found. 
Specifically, quasi-physical models are used that regard each circle as an elastic item and treat the container as a rigid hollow container~\cite{He:2013:CQO:2560985.2561717,He2015,HE201826}. If two items, or an item and the container are squeezed against each other, the whole system would have certain elastic potential energy, and by gradient descent method like BFGS we can reduce the potential energy of the system so as to remove the overlapping. Then some Basin-hopping strategy can be used to jump out of the local optimal trap where overlaps still exist.
This category mainly includes some quasi-physical  algorithms~\cite{Huang2003b,WANG2002440,Liu2016,Zhang2005,HUANG2011474}, basin hopping algorithms~\cite{Addis2008}, iterated Tabu search algorithms~\cite{Fu2012,Zeng2015}, and evolutionary search algorithms~\cite{Flores2016}. Huang \textit{et al.} proposed a global optimization algorithm based on quasi-physics, tested on instances of $1 \leq n \leq 200$ and obtained 63 better packings~\cite{HUANG2011474}.
He \textit{et al.} proposed a new quasi-physical quasi-human algorithm (QPQH)~\cite{HE201826} that utilizes the local neighbor information to speed up the calculation, tested on instances of $n=1, 2, ...320$, and obtained 66 denser layouts with smaller container radius, which is the current state-of-the-art.

 To our knowledge, there is no formal publications on instances of $n>320$, probably due to the large computational complexity. On the circle packing website \url{http://www.packomania.com}, the website maintainer Eckard Specht reported results for $n = 1$ to $5000$ for ECPP, using his “program cci, 1999–2014”. 
 However, he did not report the running time and computing machine, or release his code.


The quasi-physical model is a general model popularly used for solving ECPP, which includes a key algorithm to obtain suboptimal layout and a basin-hopping strategy to jump out of the local optimum. Our proposed method adapts this framework, and our main contribution is the design of the mini-batch BFGS method that greatly speeds up the BFGS normally used for ECPP, such that we can solve up to $n=1500$ items, and we believe this is a big progress for the general quasi-physical model. 

\vspace{-0.5em}
\section{Problem Formulation}
\label{020problem}
The equal circle packing problem (ECPP) in decision version is to ask whether we can pack $n$ unit circles into a circular container with fixed radius $R$, such that all circle items are within the border of the container and any two circle items do not overlap with each other. 

Formally speaking, we build a Cartesian coordinate system with its origin located at the center of the container and the coordinate of the center of circle $i$ is denoted by $(x_i,y_i)$, $i \in \{1,2,...,n\}$,. Then we denote any layout configuration 
 by $X=(x_1,y_1,x_2,y_2,...,x_n,y_n)$. Our purpose is to find a packing pattern of $n$ circles without overlapping, i.e., to find $(x_i,y_i)$, $i \in \{1,2,...,n\}$, such that:
\begin{equation*}
\sqrt{x_i^2 +y_i^2} + 1 \le R, ~~ 
\sqrt{(x_i-x_j)^2+(y_i-y_j)^2} \ge 2, 
\end{equation*}
where $i, j \in \{1,2,...,n\}, i\not=j$.
The first constraint 
denotes that any circle item does not intersect with the container and the second constraint
indicates that any two items do not overlap with each other. Thus, we need to find $2n$ real numbers to satisfy the two constraints, in which case we call $X$ a feasible layout. 



\vspace{-0.5em}
\section{The General Quasi-physical Model}
Among the current best approaches, researchers build a quasi-physical model to address this continuous optimization problem \cite{Huang2003b,HE201826}. Regard the container as a rigid hollow item (denoted as item ``0") fixed at the origin, and each circle $i$ as a movable elastic circular item $i$. There will be some elastic potential energy if any two elastic items overlap, 
or an item overlaps with the border of the container
Then we can calculate the elastic potential energy for a layout configuration $X$,  
and if we reduce the potential energy by some gradient descent method, there will be less overlapping among the items.

\textbf{Definition 1} \textit{Overlap Depth}. There are two kinds of overlap, circle-circle overlap and circle-container overlap. The circle-circle overlap depth is defined as:
\vspace{-0.5em}
\begin{equation}
d_{ij}=\max\left(2-\sqrt{(x_i-x_j)^2+(y_i-y_j)^2},0\right),
\end{equation}
where $i \not= j$. And the circle-container overlap depth is defined as: 
\vspace{-0.5em}
\begin{equation}
d_{0i}=\max\left(\sqrt{x_i^2+y_i^2}+1-R,0\right).
\end{equation}

\textbf{Definition 2} \textit{Elastic Potential Energy}. The elastic potential energy of the items is proportional to the square of the overlap depth. The potential energy $U_i$ of circle $i$ is defined as 
$U_i = \sum_{j=0,j\not=i}^{n}d_{ij}^2$.
And the total potential energy $U(X)$ is 
$U(X) = \sum_{i=1}^n U_i$.

Obviously, the total energy $U \ge 0$ for any layout configuration. $U=0$ if and only if $X$ is a feasible layout, i.e. $U$ is a global minimal potential.
Thus, for a fixed $R$, we minimize $U$ as the objective function so as to find a feasible solution.

\vspace{-0.5em}
\section{The Proposed SIDM Algorithm}
\label{030alg}
We adopt the general quasi-physical model for ECPP, and the key issue is how to find a local minimum of the potential energy efficiently such that we can handle large scale instances. The advantage of our method is that it can efficiently find a feasible layout, which is also a global minimum layout for a fixed container radius. 
In the following discussion, we will focus on the global optimization problem using best-known radius reported on the packomania website.

There are three procedures for a feasible layout search. First, a local search procedure finds the local minimum or saddle point, in which our stochastic item descent method is proposed. The second is the basin-hopping procedure, for which we design a flexible strategy of shrinking the container radius. Finally, the global search procedure combines the local-search and the basin-hopping procedure to search for a solution iteratively within reasonable time. 

\subsection{Stochastic Item Descent Method}
For the local search procedure, we randomly select items as a mini-batch and use the classical BFGS~\cite{liu1989limited} algorithm for gradient descent. 
The main idea of BFGS is to use the gradient information of the objective function $U$ to approximate the inverse of Hessian matrix rather than to calculate the second-order derivative at each iteration. 

For simplicity, we use $X^s$ to denote the layout of a subset of unit circles, and $U^s$ is the corresponding elastic potential energy function of this set of circles. The complementary set of $X^s$ is denoted as $X^c$. The BFGS iteration for minimizing the potential energy $U^s$ has the form:
\begin{equation}
\label{eq1}
X^s_{k+1} \leftarrow X^s_k - \alpha_k H_{k} g_{k},
\end{equation}
in which $X^s_k$ is the layout configuration at iteration $k$, $g_{k}$ is the gradient of $U^s$ at $X^s_k$, $H_{k}$ is a positive definite  approximation of $\nabla^{2}U^s(X^s_k)^{-1}$ and $\alpha_{k}$ is the step length (learning rate) at each iteration, defined in Eq. (\ref{eq2}). $H_{k}$ is updated dynamically by Eq. (\ref{eq3}), in which $I$ is the identity matrix and $u_k$, $v_{k}$ are defined in Eq. (\ref{eq4}).
\begin{equation}
\label{eq2}
\alpha_k = \mathop{\arg\min}_{\alpha \in R^+} U^s(X^s_k-\alpha H_k g_k)
\end{equation}
\vspace{-0.5em}
\begin{equation}
\label{eq4}
u_k = X^s_{k+1} - X^s_k, ~~~
v_k = g_{k+1} - g_k{}
\end{equation}
\vspace{-0.8em}
\begin{equation}
\label{eq3}
H_{k+1} = \left(I-\frac{v_k u_k^T}{u_k^T v_k}\right)^T H_k 
\left(I-\frac{v_k u_k^T}{u_k^T v_k}\right) + \frac{u_k u_k^T}{u_k^T v_k}
\end{equation}

Based on above definitions, we design a local BFGS algorithm, Algorithm \ref{alg1}, for optimizing the potential for circles in $X^s$, while other circles in $X^c$ are all fixed in the algorithm. 

\begin{algorithm}[htbp]
\caption{Local BFGS Algorithm}
\label{alg1}
\begin{algorithmic}[1]
\renewcommand{\algorithmicrequire}{\textbf{Input:}}
\renewcommand{\algorithmicensure}{\textbf{Output:}}
\REQUIRE ~~\\A layout for a subset of circles $X^s$;
           \\Container radius $R$.
\ENSURE ~~\\A local minimum layout $X^{s*}$.
\STATE iteration step $k \leftarrow 0$;
\STATE $X^s_k \leftarrow X$;
\STATE $H_k \leftarrow I$;
\STATE calculate $g_k$;
\WHILE {$k\le MaxIterNum$}
	\STATE calculate $\alpha_k$ by Eq. (\ref{eq2});
	\STATE calculate $X^s_{k+1}$ by Eq. (\ref{eq1}) ;
	\IF {$U^s \le 10^{-20}$ or $\|g_k\| \le 10^{-10}$}
		\RETURN layout $X^s_{k+1}$ as $X^{s*}$;
	\ENDIF
	\STATE calculate $g_{k+1}$;
	\STATE calculate $u_k, v_k, H_{k+1}$ by Eq. (\ref{eq4}) and (\ref{eq3}); 
	\STATE $k \leftarrow k+1$;
\ENDWHILE
\RETURN layout $X^s_k$ as $X^{s*}$.
\end{algorithmic}
\end{algorithm}

Combining the random selection of batches on unit circles with local BFGS algorithm, we have our stochastic item descent algorithm (SIDM). 
The specific idea is to randomly select a subset of circles at each time, and call local BFGS on this subset to get a locally better layout. Then we continue to randomly select another batch of circles in the remaining set and repeat such operation until all the circles have been selected in a batch. This is equivalent to a random grouping of all circles for one round of iteration, the number of circles per group is recorded as $s$ (except for the last group), and the local BFGS algorithm is called iteratively for each group.

If we continue do another random grouping on the circles at the next round of iteration and run BFGS iteratively for each group again, 
then after $k$ rounds of iterations, it is probably that the potential energy of the whole system is still relatively high. 
Therefore, we consider reducing the number of groups for each round, which means the number of circles in each group increases. The local BFGS algorithm is still applied to reduce the potential energy of each group. We need to go through $k/2$ rounds until all circles are in one group in the end,  
in which case we run the local BFGS for the whole system. As the overall packing is already relatively good, a local minimum packing layout can be quickly obtained. 
The reason why we do not choose a fixed group size but increase $s$ gradually is that small fixed group size may cause oscillation during the iterations like stochastic gradient descent for neural network training, making it hard for the potential energy to converge to a local minimum.
The pseudo code of the entire process is in Algorithm \ref{alg2}.

\begin{algorithm}[htbp]
\caption{Stochastic Item Descent Method}
\label{alg2}
\begin{algorithmic}[1]
\renewcommand{\algorithmicrequire}{\textbf{Input:}}{}
\renewcommand{\algorithmicensure}{\textbf{Output:}}
\REQUIRE ~~\\A layout configuration $X$; 
           \\Container radius $R$.
\ENSURE ~~\\A local minimum layout $X^*$.
\STATE $s \leftarrow 100$;
\STATE $k \leftarrow 10$;
\STATE $g \leftarrow \lfloor \frac{n}{s} \rfloor$;
\WHILE {$g \ge 1$}
	\FOR {$i = 1$ to $k$}
		\STATE randomly select $s$ circles as a group, with a total of $g$ groups;
		\STATE run Algorithm \ref{alg1} for each group;
		\IF {$U \le 10^{-20}$}
			\RETURN current layout as the $X^*$;
		\ENDIF
	\ENDFOR
	\STATE $s \leftarrow \min(s * 2,n)$;
	\STATE $k \leftarrow \max(\lfloor \frac{k}{2} \rfloor, 1)$;
	\STATE $g \leftarrow \lfloor \frac{n}{s} \rfloor$;
\ENDWHILE
\RETURN current layout as $X^*$.
\end{algorithmic}
\end{algorithm}

\begin{figure}[htbp]
\centering
\vspace{-1em}
\includegraphics[width=0.45\textwidth]{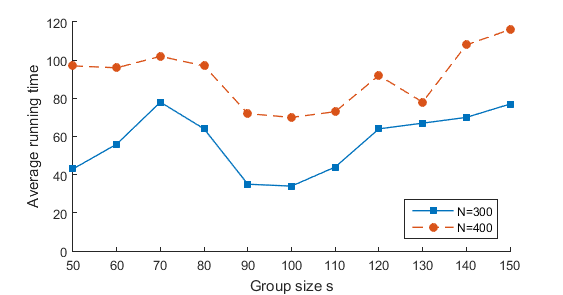}
\caption{Comparison on the average running time of 10 runs of SIDM to find a good group size.}
\vspace{-0.5em}
\label{fig2}
\end{figure}

The selection of the group size $s$ has an impact on the algorithm efficiency. We experimentally tested on two instances of $n=300$ and $n=400$ with various group sizes $s=50, 60, 70, ..., 150$. 
We compare the average running time of 10 runs that reach local minimum layout. The results are illustrated in Figure \ref{fig2}, in which we see $s=100$ is the best.

\subsection{Basin-hopping and Global Search}
The stochastic item descent usually obtains a local minimum layout or a saddle point in many cases and can not guarantee the elastic potential energy of the whole system to be small enough, aka a feasible layout may not be found. In such case, we need to consider appropriate basin-hopping strategy to help the current configuration jump out of the local optimum at the same time have a better chance to move toward the global optimum.

The shrinking strategy has a good impact on the layout with dense inner packing and sparse outer packing~\cite{HE201826}. Intuitively, if we make circles near the container center denser and make more use of the inner space, we may obtain a better layout. 
In order to get a global optimal layout, we often need to run the basin-hopping strategy multiple times. QPQH uses an identical shrinking scale for each initial shrinking radius. In practice, as the number of hops increases, it is unnecessary to squeeze the circle too far inside, and the circles near the boundary still need more precise adjustment because they are more scattered and irregular. 
Therefore, we adapt and improve the basin-hopping strategy of QPQH~\cite{HE201826} by shrinking the radius of the container more flexibly.

The coordinates of all circles are fixed and the container radius is reduced by a factor of $\gamma$ ($0 < \gamma < 1$): $R = \gamma R_0$, where $R_0$ in the initial container radius and $\gamma$ is defined as: 
\vspace{-0.5em}
\begin{equation}
\label{cal_gamma}
\gamma =\alpha + \beta \cdot hops + \frac{1 - \alpha- \beta \cdot hops}{m}k,
\end{equation}
in which $\alpha$ is the initial shrinking scale of the container radius, $hops$ is the times of running basin-hopping procedure, $\beta$ is the factor corresponding to $hops$ that adjusts the shrinking scale during the iterations, $m$ is the number of generated new layouts and $k$ varies form $0, 1, 2$ to $m-1$. Then we run stochastic item descent to reach a new layout. 

 If $\alpha$ is too small, all the circles will converge to the center of the container and most dense packing will be broken severely. If $\alpha$ is too large, there is little impact by shrinking the container radius. If $\beta$ is too small/large, the shrinking scale of each basin-hopping increases too slowly/quickly during the iteration. Besides, if $m$ is too small, the probability of generating new layouts with high quality is small; if $m$ is too large, it is very slow to generate $m$ new layouts. 
 The values are chosen empirically: $\alpha = 0.4$, $\beta = 0.03$ and $m=10$. 



    

\begin{algorithm}[htbp]
\caption{Global Search Procedure}
\label{alg4}
\begin{algorithmic}[1]
\renewcommand{\algorithmicrequire}{\textbf{Input:}}
\renewcommand{\algorithmicensure}{\textbf{Output:}}
\REQUIRE ~~\\The container radius $R_0$.
\ENSURE ~~\\A global or local minimum layout.
\STATE randomly generate an initial layout;
\STATE run SIDM to obtain an updated layout $X$;
\STATE $X^* \leftarrow X$;
\STATE $hops \leftarrow 0$;
\WHILE {$U(X^*) > 10^{-20}$ and time limit is not reached}
	\FOR {$k=0$ to $9$}
		\STATE calculate $\gamma$ by Eq. (\ref{cal_gamma});
		\STATE $R \leftarrow \gamma R_0$;
		\STATE run SIDM on layout $X^*$ with radius $R$ to generate a new layout, denoted as $X_k$;
		\STATE run SIDM on $X_k$ with radius $R_0$;
	\ENDFOR
	\IF {$\mathop{\min}_{k} U(X_k)<U(X^*)$}
		\STATE $X^* \leftarrow \mathop{\arg\min}_{k} U(X_k)$;
	\ENDIF
	\STATE $hops \leftarrow (hops + 1) \mod \lfloor\frac{1-\alpha}{\beta}\rfloor$;
\ENDWHILE
\RETURN current layout $X^*$;
\end{algorithmic}
\end{algorithm}

Combining the local search procedure with the basin-hopping procedure, we have the global search algorithm, Algorithm \ref{alg4}, that finds a feasible layout in a fixed container. It is initialized with a random layout. Then we run stochastic item descent to obtain a local minimum layout and then use the basin-hopping procedure to generate $10$ new layouts. We continue run SIDM on these packing patterns and if some packing is better than the current local minimum packing, we update the current packing. The algorithm terminates when a global minimum layout is obtained or the time limit is reached.
To show the key feature of the proposed method, we still denote the overall algorithm as SIDM. 

\subsection{Complexity Analysis}
This subsection compares the time complexity and space complexity of BFGS algorithm and local BFGS algorithm. For a fair comparison, we consider the complexity for one iteration of BFGS that all circle items update their positions once, and $m$ iterations of local BFGS where $m$ is the batch size so all the circle positions are also updated once.

Each iteration of BFGS algorithm (Simply regard $X^s$ as the layout of all circles) calculates the step length $\alpha$ by Eq. (\ref{eq2}), new layout by (\ref{eq1}), new gradient and Hessian matrix by Eq. (\ref{eq1}), and the time complexities are $O(nlog(\frac{len}{\epsilon}))$, $O(n^2)$, $O(n)$ and $O(n^2)$, respectively. Here $len$ is the length of real number interval in the line search, $\epsilon$ is the searching precision, and $nlog(\frac{len}{\epsilon})$ is the time complexity of the line search algorithm. Thus, the total time complexity is $O(nlog(\frac{len}{\epsilon}) + n^2)$. The memory mainly used by BFGS algorithm is to store the Hessian matrix, thus the space complexity is $O(n^2)$. 

The time complexity and space complexity for each batch of local BFGS algorithm are similar to BFGS algorithm, which are $O(\frac{n}{m}log(\frac{len}{\epsilon}) + (\frac{n}{m})^2)$ 
and $O((\frac{n}{m})^2)$. So for $m$ batches of local BFGS,
the time complexity is $m$ times of the complexity of a single batch of local BFGS algorithm, i.e., $O(nlog(\frac{len}{\epsilon}) + \frac{n^2}{m})$, and the space complexity is $O(\frac{n^2}{m})$. 

The time complexity of BFGS and $m$ batches of local BFGS is mainly decided by the second term, which are $O(n^2)$ and $O(\frac{n^2}{m})$, respectively. 
The time complexity of BFGS is $m$ times of the $m$ batches of local BFGS. And obviously, the space complexity of BFGS is also $m$ times of the $m$ batches of local BFGS. 

Therefore, We can conclude that SIDM using local BFGS search is more efficient than BFGS search from a complexity analysis point of view.

\vspace{-0.5em}
\section{Experimental Results}
\label{040exp}
We present our results on instances of $n=100, 200,300,$ $... ,1500$. The best-known packing results are maintained on the packomania website, where most results of ECPP are reported for $n\leq 200$ in the literature. The packomania website maintainer, Eckard Specht, also provide results using his program \texttt{cci} for $n = 1$  to 5000. But unfortunately he did not provide running time, computing machine, or code. 
To our knowledge, no result has been formally published in the literature for $n > 320$ due to the exponentially growing of computational complexity. 
The current state-of-art results formally published in the literature are from QPQH~\cite{HE201826}, which is not updated on packomania. 
Thus, we compare with QPQH~\cite{HE201826} to demonstrate the efficiency of SIDM.
 
\subsection{Experimental Setup}
SIDM is programmed in C++ programming language and implemented in Visual Studio 2017 IDE. All experiments are carried out using a personal computer with 2.5GHz CPU and 8GB RAM. 
Table \ref{tbl1} lists the key parameters of SIDM.

\begin{table}
\centering
\caption{Key parameters of the SIDM algorithm.}
\label{tbl1}
\vspace{-0.5em}
\begin{tabular}{lll}
\hline
Parameter & Description & Value \\
\hline
$s$ & Initial group size & $100$ \\
$\alpha$ & Initial shrinking factor & $0.4$ \\
$\beta$ & Shrinking scale growing factor & $0.03$ \\
$m$	& Number of new layouts & $10$ \\
\hline
\vspace{-1em}
\end{tabular}
\end{table}

\subsection{Computational Results}
Our purpose is to evaluate whether SIDM can find a global minimum layout efficiently using the reported container radius on packomania as the fixed container radius. 

\begin{figure}[h]
\vspace{-1em}
\centering
\includegraphics[width=0.45\textwidth]{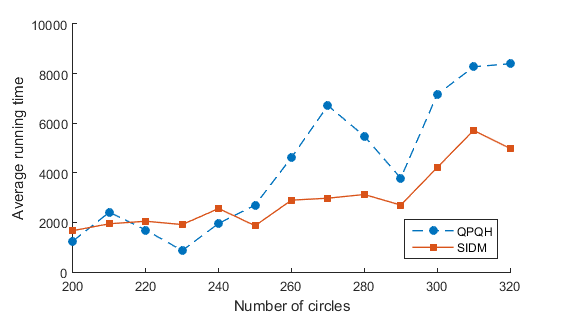}
\caption{Comparison on average running time of QPQH and SIDM.}
\vspace{-0.2em}
\label{fig5}
\end{figure}

We first compare results on instances of $n=200, 210,...,320$ between SIDM and QPQH (we use the version that the container radius is fixed). We run both algorithms for five times respectively, and show the average running time of reaching a feasible pattern in Table \ref{tbl2}. 
We also show the comparison in Figure \ref{fig5} to have an intuitive observation. 
The average running time of the two algorithms is close when the number of circles is small in 200 to 250. 
But as the number of circles increases, SIDM behaves more efficiently than QPQH. 

Then, for 15 instances of $n=100, 200, ..., 1500$, we randomly place $n$ circles in the container and run the overall SIDM algorithm. 
We will stop the search when a global minimum layout is found, or the maximum time limit of 15 hours is reached. For each instance, we run SIDM 10 times to reduce the impact of randomness.
The results listed in Table \ref{tbl3} show that SIDM can find the global minimum layout except for $n=1400$. 
The hit count indicates the number of successful times for 10 times of running, and the time indicates the average running time for successful runs.

\begin{table}
\centering
\caption{Comparison on average running time.}
\vspace{-0.5em}
\label{tbl2}
\begin{tabular}{cccc}
\hline
$n$ & $R_0$ & QPQH (s) & SIDM (s) \\
\hline
200&	15.4632748785&	\textbf{1250}&	1668	\\
210&	15.8792012772&	2412&	\textbf{1945}	\\
220&	16.2253735494&	\textbf{1690}&	2047	\\
230&	16.5964300724&	\textbf{865}	&	1912	\\
240&	16.8971658948&	\textbf{1960}&	2560	\\
250&	17.2629622393&	2697&	\textbf{1867}	\\
260&	17.6049551932&	4617&	\textbf{2897}	\\	
270&	17.8872656677&	6712&	\textbf{2976}	\\
280&	18.2472267427&	5478&	\textbf{3125}	\\	
290&	18.5493750704&	3782&	\textbf{2698}	\\	
300&	18.8135833638&	7153&	\textbf{4211}	\\
310&	19.1848594632&	8274&	\textbf{5712}	\\
320&	19.4562307640&	8397&	\textbf{4987}	\\
\hline
\end{tabular}
\vspace{-0.5em}
\end{table}

\begin{table}
\centering
\caption{Experimental results for $n=100,200,...,1500$.}
\vspace{-0.5em}
\label{tbl3}
\begin{tabular}{cccr}
\hline
$n$ & $R_0$ & Hit count & Time (s) \\
\hline
100 &	11.0821497243 &	1/10	&	2562	\\
200 &	15.4632748785 &	8/10	&	1772	\\
300	& 	18.8135833638	&	7/10	&	4326 	\\
400	&	21.6895717951	&	7/10	& 	7921 	\\
500 & 	24.1329376240	&	6/10	&	9865	\\
600	&	26.4274162694	&	4/10	&	16372	\\
700	&	28.4958443164	&	5/10	&	12369	\\
800	&	30.4212133790	&	3/10	&	15893	\\
900	&	32.2330843545	&	1/10	&	13715	\\
1000&	33.9571409147	&	1/10	&	21735	\\
1100&	35.6161932968	&	2/10	&	19816	\\
1200&	37.1121608416	&	1/10	&	34682	\\
1300&	38.6047666608	&	2/10	&	28871	\\
1400&	40.0604065845 &	0/10 &	------	\\
1500&	41.4126836805 &   1/10 &	41286	\\
\hline
\end{tabular}
\vspace{-0.5em}
\end{table}

The experimental results indicate that with the increase on number of circles, in most cases SIDM can find a feasible layout, 
and the running time increases almost linearly (2562 for $n=100$, $2562\cdot 15 = 38430$, 41286 for $n=1500$). 
By comparison, QPQH can not output any feasible results for $n = 400, 500, ...,1500$ within the time limit.

\vspace{-0.5em}
\section{Conclusion}
\label{050con}
Inspired by the idea of SGD in the area of machine learning, we propose a stochastic item descent method for large-scale equal circle packing problem (ECPP), which randomly divides the circles into batches and runs BFGS on the corresponding potential energy function in iterations. In order to obtain a solution with high quality, we increase the batch size during the iterations. Besides, we improve the basin-hopping strategy and shrink the radius of the container more flexibly. Experiment has demonstrated that the proposed method is efficient for large-scale equal circle packing problem. 


In future work, we will adapt SIDM via binary search for its optimization version problem of  minimizing the container radius, and try the SIDM idea on various circle packing problems, such as equal or unequal circles packing with various container shape. 
We also believe SIDM can be adapted for other classic optimization problems where gradient descent method has been used for optimization, including those problems occurring in the optimization process of large scale machine learning.




\bibliographystyle{named}
\bibliography{ref}


\end{document}